\input amstex
\magnification=1200
\input amsppt.sty
\loadmsbm
\NoBlackBoxes
\baselineskip 20pt
\parindent 20pt
\define \CP{\Bbb C\Bbb P}

\define \CPt{\Bbb C\Bbb P^2}
\define \Dl{\Delta}
\define \ri{\rightarrow}

\define\bnb{\bigskip \noindent \bf}
\define \Int{\operatorname{Int}}
\define \p{\partial}
\define \G{\Gamma}
\define \tm{\times}
\define \sk{\smallskip}

\define \vp{\varphi}

\define\C{{\Bbb C}}

\define\Z{{\Bbb Z}}

\define\MM{{\Cal M}}
\define\M{{\Cal M}}
\def\boxit{\vbox{\hrule\hbox{\vrule\kern5pt
             \vbox{\kern5pt}\vrule}\hrule}}

\def\today{\ifcase \month
\or January \or February \or March \or April \or May \or June \or
July \or August \or September \or October \or November \or
December\fi \thinspace
\number \day, \space \number \year}

\topmatter

\centerline{\bf Applications of braid group techniques to the decomposition}
\centerline{\bf  of moduli spaces, new examples}

\medskip

\centerline{A. Robb, M. Teicher}
\medskip
\centerline{Department of Mathematics and Computer Science}
\centerline{Bar Ilan University}
\centerline{52900 Ramat Gan, Israel}

\medskip
\centerline{robb\@bimacs.cs.biu.ac.il, teicher\@bimacs.cs.biu.ac.il}

\abstract\nofrills{\bf Abstract.}\ Every smooth minimal complex algebraic
surface of
general type, $X$, may be mapped into a moduli space, $\MM_{c_1^2(X),
c_2(X)}$, of minimal surfaces of general type, all of which
have the same Chern numbers.
Using the braid group and braid monodromy,
we construct infinitely many new examples
of pairs of
 minimal surfaces of general type which have the same Chern numbers and
non-isomorphic fundamental groups.
Unlike previous examples, our results
include $X$ for which $|\pi_1(X)|$ is arbitrarily large.
Moreover, the surfaces are of positive signature.
This supports our goal of using the braid group and
fundamental groups   to decompose
$\MM_{c_1^2(X),c_2(X)}$
into connected components.  \endabstract
\thanks This research was partially supported by the Emmy Noether Mathematics
Research Institute, Bar-Ilan University, Israel.\endthanks
\endtopmatter
\bnb 0. Introduction \rm

 It was proven by Gieseker that there exists a
quasi-projective coarse moduli space, $\MM$, of minimal surfaces of
general type.  This space is a union of components,
$\MM_{c_1^2(X),c_2(X)}$,
in which all
members have the same Chern numbers.  A major problem in the theory
of surfaces is the search for discrete invariants which characterize
the connected components of $\MM_{c_1^2,c_2}$ [Ca2], [P].

 The Chern numbers are, of course, topological invariants; $c_2(X)$
is the topological Euler characteristic, and $c_1^2(X)$ may be computed from
$c_2(X)$ and the signature of $X$.  Hence, this
may also be regarded as a problem in four-dimensional topology.

 A discrete invariant which has already been used
successfully to distinguish connected components of $\MM_{c_1^2,c_2}$
is the divisibility index, $r(X)$, of a surface, $X$.  If $K_X$ is the
canonical divisor of $X$, then $r(X)$ is the largest positive integer
such that $K_X$ is linearly equivalent to $rD$ for some divisor
$D$ of $X$.  The divisibility
index is a deformation invariant.  Catanese and Manetti have
each produced examples of minimal surfaces of general type with the same
Chern numbers and different divisibilities [Ca1].
Later Catanese and Manetti produced examples of homeomorphic minimal
surfaces of general type with the same divisibility which are not
deformations of each other, i.e., they are in different connected components
of moduli spaces (see [Ca3], [Ma1], [Ma2]).  All of their surfaces are
simply-connected.

 Another discrete invariant is the fundamental group.
There exist pairs of surfaces which have the same Chern
numbers and non-isomorphic fundamental groups.  For example, the
Godeaux construction can be used to produce surfaces, $X$, such that
$c_1^2(X) = 2$, $c_2(X) = 10$, and $\pi_1(X)$ is one of the following:
$\oplus_3 \Z_2$; $\Z_2 \oplus\Z_4$; $\Z_8$; the multiplicative group
$\{ \pm 1, \pm i \pm j \pm k\} \subset \Bbb H$ [BPV].  However, there
are few examples of such pairs of surfaces and the possible fundamental
groups are all relatively small.

\medskip In this article, we use new results on Galois covers of
Hirzebruch surfaces [MoRT] to construct infinitely many new examples of
pairs of minimal surfaces of general type, $X$ and $Y$, such that $$c_1^2(X) =
c_1^2(Y) \qquad \qquad c_2(X) = c_2(Y) \qquad \qquad \pi_1(X) \not \cong
\pi_1(Y).$$

 In all of our examples, $X$ and $Y$ have positive signature and $$\pi_1(Y)
= \{0\}.$$
The fundamental group of $X$ is finite but may be arbitrarily large.
Indeed,  for every $n \in \Z$,
there exists a pair $X$, $Y$, such that $$
|\pi_1(X)| > n.$$

This article is organized as follows.  Section 1 is an introduction to the
braid monodromy map associated to an algebraic curve.
Section 2 describes the Galois cover of an algebraic surface (our examples
are Galois covers of Hirzebruch surfaces).
Section 3 describes how to use braid monodromy theory to calculate the
fundamental group of a Galois cover.
Section 4 contains our examples of pairs of surfaces with the same Chern
numbers and non-isomorphic fundamental groups.
Section 5 gives more examples, which are based on  a work in
preparation by the authors.
It also describes how more examples might be obtained using Galois covers
of K3 surfaces.

 \bigskip \subheading{1. The Braid Group Techniques}

We introduce braid monodromy, which is used to compute fundamental groups
related to surfaces of general type:  fundamental groups of the complement
of a branch curve and fundamental groups of Galois covers.

We consider the following situation.

$S$ is a curve in $\C^2,$ \ $p=\deg S.$

$\pi: \C^2\ri \C$ a projection on the first coordinate.

$N=\{x\bigm| \#\pi^{-1}(x)\cap S\lvertneqq p\}.$

 Let $E$   be a closed disk on $x$-axis   s.t.
 $N\subset\Int(E).$  We choose $u\in\p E.$

Let $D$ be a closed disk on the $y$-axis s.t. $\pi^{-1}(E)\cap S\subset E\times
D.$

$K=\{y\bigm| (u,y)\in S\} \ =\{q_1,\dots,q_p\}.$

In such a situation, we are going to introduce ``braid monodromy."
\smallskip
\definition{Definition} \ \underbar{Braid monodromy of $S$ w.r.t. $E\times
D,\pi, u.$}

Every loop  in $E-N$ starting at $u$ has liftings to a system of $p$
paths in $(E-N)\times D$ starting at $q_1,\dots, q_p.$
Projecting them to $D$ we get $p$ paths in $D$ defining a motion
$\{q_1(t),\dots,q_p(t)\}$ of $p$ points in $D$ starting and ending at $K.$

This motion defines a braid in $B_p[D,K],$ as explained in Chapter III of
[MoTe4]. Thus we get a map $\vp:\pi_1 (E-N,u) \ri B_p[D,K].$
This map is evidently a group homomorphism, and it is the braid
monodromy of $S$ w.r.t. $E\tm D, \pi, u.$
We sometimes denote $\vp$ by $\vp_u.$\enddefinition
\smallskip
\definition{Definition} \ $\underline{\text{Braid monodromy of}\  S \
\text{w.r.t.} \ \pi,u.}$

When considering the braid induced from the previous motion as an element of
the group $B_p [\C_u, K]$ we get the homomorphism $\vp:\pi_1(E-N,
u)\ri B_p[\C_u, K]$ which  is called the braid mondromy of $S$ w.r.t. $\pi,
u.$
We sometimes denote $\vp$ by $\vp_u.$\enddefinition
\sk

\proclaim{Proposition - Example 1.1} \ Let $E=\{x\in\C ||x|\leq
1\},$\ $D=\{y\in\C |y \leq R\},$ $R\gg 1, S$ is the curve $y^2=x^\nu,\ u=1.$
Clearly, here $n=2, N=\{0\}, K=\{-1, +1\}$ and $\pi_1(E-N, 1)$ is generated by
$\G=\p E$ (positive orientation).
Denote by \linebreak $\vp:\pi_1(E-N, 1)\ri B_2[D,K]$ the braid monodromy of
$S$ w.r.t. $E\times D, \pi, u.$

Then $\vp(\G)=H^\nu,$ where $H$ is the positive half-twist defined by $[-1,
1]$\linebreak (``positive generator" of $B_2[D, K]$).\endproclaim
\smallskip
\demo{Proof} \ We can write $\G = \{e^{2\pi it}, t\in[0,1]\}.$
Lifting $\G$ to $S$ we get two paths:
$$\align
\delta_1(t) &= \left(e^{2\pi it}, \ e^{2\pi i\nu t/2}\right)\\
\delta_2(t) &= \left(e^{2\pi it}, \ -e^{2\pi i\nu t/2}\right).\endalign$$

Projecting $\delta_1(t),$ $\delta_2(t)$ to $D$ we get two paths:
$$\alignat 2
&a_1(t) =e^{\pi it\cdot \nu}, \qquad \qquad  &&0\leq t\leq 1\\
&a_2(t) = -e^{\pi it\cdot \nu}, &&0\leq t\leq 1.
\endalignat$$

This gives a motion of \{1, -1\} in $D.$
This motion is the $\nu$-th power of the motion $\Cal M$:
$$\alignat 2
&b_1(t) =e^{\pi it}, \qquad \qquad  &&0\leq t\leq 1\\
&b_2(t) = -e^{\pi it}, &&0\leq t\leq 1.
\endalignat$$

The braid of $B_2\left[D, \{1, -1\}\right]$ induced by ${\M}$ coincides with
the half-twist $H$ corresponding to $[-1, 1]\subset D.$ Thus $\vp(\G)=H^\nu.$
\quad \qed \enddemo
\sk
\proclaim{Proposition - Example 1.2} \ Let $S$ be a union of
$p$ lines, meeting in one point $s_0, s_0=(x(s_0),y(s_0)).$
Let $D, E, u, K$ be as before.
Let $\vp$ be the braid monodromy of $S$ w.r.t. $ E\times D, \pi,u.$
Clearly, here $N=$ single point $x(s_0)$ and $\pi_1(E-N, u)$ is generated by
$\G=\p E.$
Then $\vp(\G)=\Dl^2_p=\Dl^2_p \left[ D, K(u)\right].$\endproclaim
\smallskip
\demo{Proof} \ By a continuous change of $s_0$ and the $n$ lines passing
through
$s_0$ (and by uniqueness of $\Dl^2_p$) we can reduce the proof to the following
case: $S=UL_k,$ \linebreak $L_k\:$ $y=j_kx, \ j_k=e^{2\pi ik/p}$ \quad
$k=0\dots p-1.$ Then $N=\{0\}.$
We can take $E=\{c| |x|\leq 1\},$ $u=1,$\ $\G=\p E=\{x=e^{2\pi it},$
$t\in[0,1]\}.$
Lifting $\p E$ to $S$ and then project it to $D$ we get $n$ loops:
$$a_k(t)=e^{2\pi i(t+k/p)}\quad k=0\dots p-1\quad t\in[0,1].$$

Thus the motion of $a_k(0)$ represented by $a_k(t)$ is a full twist which
defines the braid $\Dl^2_p\left[ D, \{a_k(0)\}\right] =
\Dl^2_p\left[ D, K(1)\right].$

(To check the last fact, see the corresponding actions in
$\pi_1\left(D-K,u\right)$). \ \qed\enddemo

Let $S$ be a   curve in $\Bbb C^2,$\ $p=\operatorname{deg} S,$
\ $\Bbb C_u=\{(u,y)\}.$
There exists an epimorphism $\pi_1(\Bbb C_u-S,u_0)\to \pi_1(\Bbb
C^2-S,u_0),$ so a  set of generators for $\pi_1(\Bbb C_u-S,u_0)$ determines a
set of generators for $\pi_1(\Bbb C^2-S,u_0).$

There is a classical theorem of Van Kampen from the 30's, ([VK], [Z]), which
states that all relations in $\pi_1(\Bbb C^2-S,u_0)$ come from the braid
group $B_p$ via the braid monodromy $\varphi_u$ of $S.$
We shall formulate it precisely.

Choose $L$, a line in infinity transverse to $S.$
Let $\Bbb C^2=\Bbb C\Bbb
P^2-L.$

Choose coordinates $x,y$ on $\Bbb C^2.$
Let $\varphi_u:\pi_1(\Bbb C-N,u)\to B_p$ the braid monodromy of $S$ with
respect to $\pi,u.$

The group $\pi_1(\Bbb C_u-S,u_0)$ is a free group.

\proclaim{Van Kampen Theorem}

 $\pi_1(\Bbb C^2-S,u_0)\simeq \pi_1(\Bbb
C_u-S,u_0) /\{\beta(V)=V\bigm|\beta\in Im\ \varphi_u,V\in\pi_1(\Bbb
C_u-S,u_0)\}. $ \endproclaim

The above formulation of Van Kampen is not very practical because the group
presentation is not finite.  It is possible to simplify the presentation so
that it is finite.  See, for example, [L]. Moreover, since we consider
branch curves
which are cuspidal, one can formulate Van Kampen's theorem with relations of
types $AB,$\ $ABA^{-1}B^{-1}$ and $ABAB^{-1}A^{-1}B^{-1}.$ Even this
presentation is rather long and complicated, and in order to apply it we have
to find symmetries in the braid monodromy factorizations. \bigskip

\subheading{2. Galois Covers and their Chern Numbers}

 We use the Galois cover construction of Miyaoka to construct our examples.

\definition {Definition}\ $\underline{\text{Galois cover}}$

For $X$ a surface in $\Bbb C\Bbb P^N$ and $f:X\to \Bbb C\Bbb P^2$ a generic
projection, we define the Galois cover of $X$ and $f$ w.r.t.the full
symmetric group  as: $\tilde
X=X_{\operatorname{Gal}}=\overline{(X\underset
\underbrace{\CPt\quad\CPt}_{n\ \text{times}}\to{\times\dots\times}
X)-\Dl}$
where $\Dl$ is the set of $n$-tuple $(a_1\dots a_n)$,\ $a_i=a_j$ for some
$i\ne j.$ Let $\tilde f: \tilde X\to\CPt$  denote the natural projection.

  Let $S \subset \CP^2$ denote
the branch curve of $f$.  The curve $S$ is singular, with ordinary
singularities - nodes and cusps.

  The surface $\tilde X$ is smooth.  If $\deg (S) > 6$, then $\tilde X$ is
minimal and of general type [MoTe2].  \enddefinition

\proclaim{2.1 Lemma} Let $n =\deg(X)$ and $m =\deg(S)$.  Let $d$ and $\rho$
denote the respective numbers of nodes and cusps of $S$.  Then
 $$c_1^2(\tilde X)
= \frac {n!} 4 (m-6)^2 \qquad \qquad c_2(\tilde X) = n! \left ( \frac 1 2
m(m-3) + 3 - \frac {3d} 4 - \frac {4\rho} 3 \right).$$\endproclaim

\demo{Proof} \rm See 7.1.1 of [MoRoTe].  \boxit
\enddemo
\proclaim{2.2 Lemma
} Let $E$ and
$K$ denote respective hyperplane and canonical divisors of
$X$.  Then the Chern numbers of $\tilde X$ are functions of
$c_1^2(X)$, $c_2(X)$, $\deg (X)$, and $E\cdot K$.\endproclaim

\demo{Proof}  Let $g$ denote the genus of an algebraic curve and let $e$
denote the topological Euler characteristic of of a space.

Let $R \subset X$ denote the ramification locus of $f$.  The
curve $R$ is a non-singular model of $S$.
By the Riemann-Hurwitz formula, $R = K
+ 3E$.   Thus
$$
e(R)  = -R \cdot (R + K)
=-(K+3E) \cdot (2K + 3E)$$
It follows that
$g(R)$ is determined by $n = E^2$, $c_1^2(X) = K^2$, and $K \cdot E$.
Similarly,  $e(E)$ is determined by these quantities.
Because $m = \deg (S) = \deg (R) = E \cdot (K
+ 3E)$, we have that $m$ is determined by these quantities.

  Let $S^*$ denote the dual curve to $S$ and let $\mu = \deg (S^*)$.
By the preceding paragraph, and by Lemma 2.1, it suffices to show that
$\mu$, $d$, and $\rho$ are determined by $m$, $n$, $e(E)$, $c_2(X)$, and
$g(R)$.  We show this by presenting three linearly independent formulae:
$$m(m-1) = \mu + 2d + 3 \rho$$
$$g(R) = \frac {(m-1)(m-2)} 2 - d - \rho$$
$$c_2(X) + n = 2e(E) + \mu$$
The first two are Classical Plucker formulae.  For the third, we may
find a Lefschetz pencil of hyperplane sections of $X$ whose union is
$X$.  Thus, $$e(X) + n = e(\CP^1) \cdot e(E) +
\text {(number of singular curves in the pencil)}$$
where $e$ is the topological Euler characteristic.  The number of singular
curves is equal to $\mu$.  \boxit\enddemo

\demo{2.3 Remark}
Lemma 2.1 can easily be modified to give explicit formulae for
$c_1^2(\tilde X)$ and $c_2(\tilde X)$ in terms of $c_1^2(X),$ $c_2(X),$ $n,$
and $E.K.$
However, such formulae are not necessary for our result.\enddemo \bigskip

\subheading{3. Fundamental groups of Galois cover}

Let $X$ be a surface in $\CP^N.$
Let $f$ be a generic projection $f: X\to \CPt.$
Let $X_{\operatorname{Gal}}$  be the Galois cover of $X$ and $f$ w.r.t. the
full symmetric group.

 Consider  the natural homomorphism $\pi_1(\Bbb
C^2-S,u_0)\overset\psi\to \rightarrow S_n$ for $S$ the branch curve of $f$
and $u_0$ any point not in $S.$
In fact, lifting a loop at $u_0$ to $n$ paths in $X $ \
$(n=\operatorname{deg} f),$ induces a permutation of $f^{-1}(u_0).$
Since $\#f^{-1}(u_0)=n$ we thus get an element of $S_n.$
Because $f$ is a generic projection,  we have $1\to \ker\psi\to\pi_1(\Bbb
C^2-S,u_0)\to S_n\to 1.$

 To obtain an isomorphic form of $\pi_1(X_{\operatorname{Gal}})$
related to the braid monodromy we have to choose a certain system of
generators for $\pi(\Bbb C^2-S,u_0).$

Let $\pi: \Bbb C^2-\Bbb C$ be the projection on the first coordinate.
Let $M'\subseteq S\subseteq \Bbb C^2$ be the points of $S$ where
$\pi\bigm|_S$ is not etale.
Let $M=\pi(M')$.  The set $M$ is finite.
Let $u\in\Bbb C_M.$
\ $(\pi\bigm|_{S})^{-1}(u)$ is a ``good'' fibre.
Let us take $u$ real, ``far enough'' from the ``bad'' points.
Let $u_0$ be a point in $\Bbb C_u=\pi^{_1}(u),\ u_0\notin S.$

Let $S\cap \Bbb C_u=\{q_1,\dots q_p\}.$
Let $\gamma_j$ be paths from $u_0$ to $q_j,$ such that the $\gamma_j$'s do
not meet each other in any point except $u_0.$
Let $\eta_j$ be a small circle around $q_j.$
Let $\gamma_j'$ be the part of $\gamma_j$ outside $\eta_j.$
Take $\Gamma_j=\gamma_j'g_j(\gamma_j')^{-1}$.
The set $\{\Gamma_j\}$ freely generates $\pi_1(\Bbb C_u-S,u_0).$
\ $\{\Gamma_j\}_{j=1}^p$ is called a good system of generators for
$\pi_1(\Bbb C_u-S,u_0).$

We have a surjection $\pi_1(\Bbb C_u-S,u_0)\overset
\nu\to\rightarrow\pi_1(\Bbb C^2-S,u_0)\to 0.$
The set $\{\nu(\Gamma_j)\}$  generates $\pi_1(\Bbb C^2-S,u_0).$
By abuse of notation, we shall denote $\nu(\Gamma_j)$ by $\Gamma_j.$

Since $f$ is stable, $\Gamma_j$ induces a transposition in $S_n.$
So $\Gamma_j^2\in\ker\psi.$
Let $\langle \Gamma_j^2\rangle$ be the normal subgroup generated by
$\Gamma_j^2.$
Then $\langle \Gamma_j^2\rangle\subseteq \ker\psi.$
By the standard isomorphism theorems, we have:
$$1\to \ker\psi/\langle\Gamma_j^2\rangle\to\pi_1(\Bbb C^2-S,u_0)/\langle
\Gamma_j^2\rangle\overset\psi\to\rightarrow S_n\to 1.$$

It is convenient to replace $\CPt$ by its ``generic'' affine part.
Let $S$ be now the branch curve of $f$ in $\Bbb C^2.$
Let $X_{\operatorname{Gal}}^{\operatorname{Aff}}$ be the part of
$X_{\operatorname{Gal}}$ lying over $\Bbb C^2(\subseteq \CPt).$
It is evident that $X_{\operatorname{Gal}}^{\operatorname{Aff}}\to
X_{\operatorname{Gal}}$ is surjective.
We consider first $\pi_1(X_{\operatorname{Gal}}^{\operatorname{Aff}}).$

\proclaim{3.1 Proposition}
$\pi_1(X_{\operatorname{Gal}}^{\operatorname{Aff}})$ is isomorphic to
$\ker\psi/\langle\Gamma_j^2\rangle.$\endproclaim
\demo{Proof} in [MoTe1].\ \qed\enddemo

Proposition 3.1 reduces the problem of computing
$\pi_1(X_{\operatorname{Gal}}^{\operatorname{Aff}})$
 to the computation of a subgroup of  $\pi_1(\Bbb
C^2-S,u_0)/\langle \Gamma_j^2\rangle.$
 Using the braid monodromy we get via
the Van Kampen method (see \S1) a finite presentation of $\pi_1(\Bbb
C^2-S,u_0)/\langle \Gamma_j^2\rangle.$
>From this presentation via the Reidermesiter-Schreier method
we get   a finite presentation of its subgroup
$\ker\psi/\langle\Gamma_j^2\rangle$  which is
$\pi_1(X_{\operatorname{Gal}}^{\operatorname{Aff}})$. Passing to the
projective case means adding one relation which is
$\prod\limits_j\Gamma_j=1.$ \bigskip

 \subheading{4. The New Examples}

Let $F_k$ denote the Hirzebruch (rational ruled) surface of order
$k$.  The Picard group of $F_k$ is generated by divisors $C$ and $E_0$,
where $C \subset F_k$ is a fiber and $E_0 \subset F_k$ is
a zero section.  We have $E_0^2 = k$, $E_0 \cdot C = 1$, and $C^2 = 0$.
Let $K_k$ denote the canonical divisor of $F_k$.  It is well known
that for any $k$,
$$K_k = -2E_0 + (k-2)C \qquad \qquad
c_1^2(F_k) = 8 \qquad \qquad c_2(F_k) = 4$$

 Let $F_{k, (a,b)}$ denote the image
of $F_k$ under the embedding induced by $aC+bE_0$\ $(a,\ b \ge 1)$.  It is
elementary that  $$\align&\text{(a)}\ \deg(F_{k, (a,b)}) = (aC + bE_0)^2 =
2ab + b^2k\\
&\text{(b)}\ K_k \cdot (aC + bE_0) =(-2E_0+(k-2)C)(aC+bE_0)= -2a-2b-bk.
\endalign$$

For every pair of positive integers $s,\ t$,
define
$$X_{s,t} = \tilde F_{0, (s+t, 2t)} \qquad \qquad
Y_{s,t} = \tilde F_{1, (s, 2t)}.$$

\proclaim{4.1 Proposition }  For any $s,\ t$,

$$c_1^2(X_{s,t}) = c_1^2(Y_{s,t}) \qquad \qquad
c_2(X_{s,t}) = c_2(Y_{s,t}).$$\endproclaim

\demo{Proof} Since $X_{s,t}$ and $Y_{s,t}$ are Galois covers we can prove
the proposition by using Lemma 2.2.  Because all rational ruled surfaces have
the same Chern numbers, then in order to apply Lemma 2.2 it suffices to prove
that the degrees and the  intersection of the hyperplane divisor with
the canonical divisor, are the same for $F_{0, (s+t, 2t)}$ and
$F_{1, (s, 2t)}$.  This is easily done with the above formulae (a) and (b).
\boxit\enddemo

\proclaim{4.2 Theorem} Let $s$, $t$, be odd integers such that $\gcd(s,t) =
1$.  Let $n(s,t) = \deg (F_{0, (s+t,2t)}) =  4st+4t^2$.  Then
$$\pi_1(X_{s,t}) \cong \bigoplus_{n(s,t) - 2} \Z/2\Z\qquad \qquad
\pi_1(Y_{s,t}) \cong \{0\}.
$$\endproclaim
\demo{Proof} \rm The case of $X_{s,t}$ follows directly from [MoTe3],
Theorem 10.2.  This states that
$\pi_1(\tilde F_{0,(a,b)}) \cong \oplus_{n(a,b)
- 2} \Z/c\Z$ where $c = \gcd(a,b)$.
The case of $Y_{s,t}$ follows directly
from [MoRoTe], Theorem 0.1.  This states that for any $k$, $\pi_1(\tilde
F_{k,(a,b)})$ is trivial when $\gcd(a,b) = 1.$  \boxit\enddemo

\proclaim{4.3 Theorem} $X_{s,t}$ and $Y_{s,t}$ are minimal surfaces of general
type which are $4$-manifolds with positive signature.\endproclaim
\demo{Proof} \cite{MoRoTe}.\enddemo

 By combining Proposition 4.1, Theorem 4.2 and Theorem 4.3, we obtain our
examples.
Each example is a pair of 2 smooth minimal surfaces of general type  with
positive signature which have the same Chern numbers and
non-isomorphic fundamental groups.
One fundamental group is trivial and the other one is a commutative finite
group.
The examples include $X$ for which $|\pi_1(X)|$ is arbitrarily large. \bigskip
\subheading{5. Other Examples}

 It is possible to use Galois covers of other embeddings of Hirzebruch
surfaces.  For example, we obtain non-isomorphic fundamental
groups if $s$ and $t$ are both odd but not necessarily relatively prime.
Non-isomorphic fundamental groups are also obtained if $s = t$.

 However, these cases require Theorem 5.1, whose proof is in [FRT]. This
theorem is a generalization of the results of [MoTe3] and [MoRoTe].  Its
proof uses the topological and group-theoretic techniques of these articles.

\proclaim{5.1 Theorem} Let $\tilde F_{k,(a,b)}$ denote the Galois cover of
$F_{k,(a,b)}$.  Let $n(a,b) = \deg(F_{k,(a,b)})$ and $c = \gcd(a,b)$
then $$\pi_1(\tilde F_{k,(a,b)}) \cong \bigoplus_{n(a,b) - 2}
\Z/c\Z.$$\endproclaim

 It is likely that the Galois construction can be used to
produce many more examples of pairs of surfaces with the same Chern
numbers and non-isomorphic fundamental groups.  To illustrate this, we will
consider covers of K3 surfaces.
\medskip An analog of Theorem 5.1 is Conjecture 5.2.
\proclaim{5.2 Conjecture} Let $X \subset \CP^N$ be an embedded \rom{K3}
surface, with
hyperplane section $E$ and degree $n$.   Let $D \subset X$ be the
hyperplane
section of a embedding of $X$ of minimal degree, and assume that $E
\equiv sD$.
Then
$$\pi_1(\tilde X) \cong \bigoplus_{n-2} \Z_s$$.\endproclaim

Conjecture 5.2 can be proven in the
cases for which $D^2 = 4,6,8$ (in these cases,
$X$ is a complete intersection).  It
should be possible to combine results on the degenerations of K3 surfaces
of
Ciloberto, Lopez, and Miranda [CiLMi] with  braid monodromy techniques
([MoTe1], [MoTe3], [MoRoTe) to prove
the other cases.  It should be mentioned that similar results exist for
other
surfaces [R2].

  Let $a$,
$b\in \Z^+$ be distinct even square integers.  Let $X$ and $Y$ be
embedded K3
surfaces, each of degree $a^2b^2$.  Assume that the smallest possible
degrees of embeddings of $X$ and $Y$ are $a^2$ and $b^2$ respectively. It
is
elementary that the K3 surfaces $X$ and $Y$ above
satisfy the hypotheses of Proposition 4.1.  Hence, $\tilde X$ and
$\tilde Y$ have the same Chern numbers.  Conjecture 5.2 would imply that
$\pi_1(\tilde X) \cong \oplus_{n-2} \Z_b$ and
$\pi_1(\tilde X) \cong \oplus_{n-2} \Z_a$, which
would imply the existence of
many more families of examples.

\end